\documentclass[a4, 12pt]{amsart}
\usepackage{fullpage}

\usepackage[title]{appendix}
\usepackage{graphicx}
\usepackage{subcaption}
\usepackage{amsthm}
\usepackage{mathrsfs} 
\usepackage{amsfonts}
\usepackage{amssymb}
\usepackage{diagbox}
\usepackage{mathtools}
\usepackage{longtable} 
\usepackage{MnSymbol}
\usepackage{indentfirst}
\usepackage[all]{xy}
\usepackage[colorlinks=true,linkcolor=blue]{hyperref}
\usepackage{ytableau}
\usepackage{tikz-cd}
\usepackage{enumerate}
\theoremstyle{plain}
\newtheorem{theorem}{Theorem}[section]
\newtheorem{thm}{Theorem}
\newtheorem{prop}{Proposition}
\newtheorem{conjecture}{Conjecture}[section]
\newtheorem{problem}{Problem}[section]

\newtheorem{lemma}[theorem]{Lemma}
\newtheorem{corollary}[theorem]{Corollary}
\newtheorem{proposition}[theorem]{Proposition}
\newtheorem{definition}{Definition}
\theoremstyle{remark}

\newtheorem{remark}{Remark}[section]

\makeatletter

\newcommand{\Rmnum}[1]{\expandafter\@slowromancap\romannumeral #1@}
\makeatother

\def\rd{{\rm d}}

\newcommand{\qbm}[2]{\left[ {#1 \atop #2} \right]}
\newcommand{\qbml}[2]{{ {#1 \atop #2} }}

\def\bp{\mathbb{P}}

\def\br{\mathbb R}
\def\bn{\mathbb N}
\def\bz{\mathbb Z}
\def\bc{{\mathbb C}}

\newmuskip\pFqmuskip

\numberwithin{equation}{section}
\allowdisplaybreaks[2] 

\begin{document}

\title[Unimodality of bivariate formal Laurent series]{Unimodality and certain bivariate formal \\ Laurent series}
\author{Nian Hong Zhou}

\address{School of Mathematics and Statistics, The Center for Applied Mathematics of Guangxi,
Guangxi Normal University, Guilin, 541006, Guangxi, PR China}
\email{nianhongzhou@outlook.com; nianhongzhou@gxnu.edu.cn}%

\thanks{This paper was partially supported by the National Natural Science Foundation of China (No. 12301423), and the Key Laboratory of Mathematical Model and Application (Guangxi Normal
University), Education Department of Guangxi Zhuang Autonomous Region.}%
\subjclass{Primary 05A20;  Secondary  05A15, 05A17}%
\keywords{Unimodality; Gauss polynomials; refined partitions; ranks; Betti numbers}%

\begin{abstract}
In this paper, we examine the unimodality and strict unimodality of certain formal bivariate Laurent series with non-negative coefficients. We show that the sets of these formal bivariate Laurent series form commutative semirings under the operations of addition and multiplication of formal Laurent series. This result is used to establish the unimodality of sequences involving Gauss polynomials and certain refined color partitions. In particular, we solve an open problem posed by Andrews on the unimodality of generalized Gauss polynomials, establish an unimodal result for a statistic of plane partitions, and establish many unimodal results for rank statistics in partition theory, including the rank statistics of concave and convex compositions studied by Andrews, as well as certain unimodal sequences studied by Kim-Lim-Lovejoy. Additionally, we establish the unimodality of the Betti numbers and Gromov-Witten invariants of certain Hilbert schemes of points.
\end{abstract}
\maketitle


\section{Introduction}
A \emph{(strictly) unimodal sequence} is a finite sequence of real numbers that first (strictly) increases and then (strictly) decreases.  A polynomial $f(x)=\sum_{0\le i\le n}c_i x^i$ is called a symmetric (strictly) unimodal polynomial if the sequence $(c_i)_{0\le i\le n}$ is a (strictly) unimodal sequence and $c_i=c_{n-i}$ for all $0\le i\le n$. Symmetric unimodal polynomials frequently appear in combinatorics, algebra and algebraic geometry, and have been the subject of extensive research in recent years.  We refer the reader to Stanley \cite{MR1110850} and Br\"{a}nd\'{e}n \cite{MR3409348} for excellent surveys of many of these techniques and further references.

\medskip

It is well-known that the product of two symmetric unimodal polynomials with non-negative real coefficients is also a symmetric unimodal polynomial with non-negative real coefficients. This means that there exists a \emph{nice algebraic structure} for these polynomials. To understand this, we consider a special case of Laurent polynomials over $\mathbb{R}$ as defined in Definition \ref{defm} below.

\begin{definition}\label{defm}
For $\nu\in\{1,2\}$, define $\mathscr{U}_x^\nu$ as the set of all Laurent polynomials $f(x)=\sum_{i}c_i x^i$ such that $c_{-i}=c_i$ for all integers $i\ge 0$ and
\begin{equation*}
\begin{cases}c_{r+\ell\nu}\ge c_{r+(\ell+1)\nu}\;\; &0\le \ell<s_r\\
\quad c_{r+\ell\nu}=0\;\;  &\quad\ell>s_r,
\end{cases}
\end{equation*}
for each integer $0\le r<\nu$ and some integer $s_r\ge 1$. If all the inequalities ``$c_{r+\ell\nu}\ge c_{r+(\ell+1)\nu}$" are replaced with ``$c_{r+\ell\nu}>c_{r+(\ell+1)\nu}$", we refer to the set of such polynomials as $\mathscr{T}_x^\nu$.
\end{definition}
We observe that the polynomials described in Definition \ref{defm} are essentially an extended and rewritten version of symmetric (strictly) unimodal polynomials with nonnegative coefficients. Let $\nu\in\{1,2\}$, then it is clear that $\mathscr{U}_x^1\subseteq \mathscr{U}_x^2$, $\mathscr{T}_x^1\subseteq \mathscr{T}_x^2$ and $\mathscr{T}_x^\nu\subset \mathscr{U}_x^\nu$. Moreover, if $f\in \mathscr{U}_x^\nu$ (or $f\in \mathscr{T}_x^\nu$) then the sequence $(c_{\nu n+j})_{n\in \bz}$ is unimodal (or strictly unimodal) for each integer $0\le j<\nu$. 
Recall that a \emph{commutative semiring} is an algebraic structure satisfying the usual axioms for a commutative ring without requiring that addition be invertible. Then, we have the following proposition.
\begin{prop}\label{main0}
For each $\nu \in \{ 1, 2 \}$, the sets $\mathscr{U}_x^\nu$ and $\mathscr{T}_x^\nu$ both form commutative semirings under the operations of addition and multiplication of Laurent polynomials.
\end{prop}
\begin{proof}
Let $\nu\in\{1,2\}$. By definition, it is clear that the sets $\mathscr{U}_x^\nu$ and $\mathscr{T}_x^\nu$ are both closed under the addition of Laurent polynomials. It remains to show that the sets $\mathscr{U}_x^\nu$ and $\mathscr{T}_x^\nu$ both are also closed under the multiplication of Laurent polynomials. To do this, we will use similar arguments to those of Andrews \cite[Theorem 3.9]{MR1634067}. For any
$$f(x)=\sum_{n}a_nx^n\in \mathscr{U}_x^\nu\;\;\text{and}\;\; g(x)=\sum_{n}b_nx^n\in\mathscr{U}_x^\nu,$$
write $f(x)g(x)=\sum_{n}c_nx^n$. Notice that $a_{r}=a_{-r}$ and $b_r=b_{-r}$ for any $r\ge 0$, then we have
\begin{align*}
c_n=\sum_{r\in\bz}a_rb_{n-r}=\sum_{r\in\bz}a_{-r}b_{-n+r}=\sum_{r\in\bz}a_{r}b_{-n-r}=c_{-n}.
\end{align*}
Moreover, for any $n\ge 0$ we have
\begin{align*}
c_n-c_{n+\nu}=&\sum_{r\in\bz}a_r\left(b_{n-r}-b_{n+\nu-r}\right)\\
=&\sum_{r\in\bz}a_{n-r}\left(b_{r}-b_{r+\nu}\right)\\
=&\sum_{r\ge 0}a_{n-r}\left(b_{r}-b_{r+\nu}\right)+\sum_{r>0}a_{n+r}\left(b_{r}-b_{r-\nu}\right)\\
=&\sum_{r\ge 0}a_{n-r}\left(b_{r}-b_{r+\nu}\right)+\sum_{r>-\nu}a_{n+r+\nu}\left(b_{r+\nu}-b_{r}\right)\\
=&\sum_{r\ge 0}\left(a_{n-r}-a_{n+r+\nu}\right)\left(b_{r}-b_{r+\nu}\right)+\sum_{1\le r<\nu}a_{n-r+\nu}\left(b_{\nu-r}-b_{r}\right).
\end{align*}
Since $\nu\in\{1,2\}$, we have $\sum_{1\le r<\nu}a_{n-r+\nu}\left(b_{\nu-r}-b_{r}\right)=0$. Thus
\begin{align}\label{eqmmmm}
c_n-c_{n+\nu}&=\sum_{r\ge 0}\left(a_{n-r}-a_{n+r+\nu}\right)\left(b_{r}-b_{r+\nu}\right)=\sum_{t\ge 0}\left(b_{n-t}-b_{n-t+\nu}\right)\left(a_{t}-a_{t+\nu}\right).
\end{align}
Since $\nu\le n+r+\nu-|n-r|\equiv 0\pmod \nu$ for $\nu\in\{1,2\}$ and $r\ge 0$, we have
$$a_{n-r}-a_{n+r+\nu}=a_{|n-r|}-a_{|n-r|+(n+r+\nu-|n-r|)}\ge a_{|n-r|}-a_{|n-r|+\nu}.$$
This implies that $c_n\ge c_{n+\nu}$ for any $n\ge 0$ because of $b_{r}\ge b_{r+\nu}$ and $a_{|n-r|}\ge a_{|n-r|+\nu}$ for any $r\ge 0$. Additionally, if $f, g\in \mathscr{T}_x^\nu$ and $c_n>0$, then there exists $r\in\bz$ such that
$$a_{n-r}b_r=a_{t}b_{n-t}>0,$$
where $t=n-r$. In view of \eqref{eqmmmm}, if $r<0$ then $t=n-r>0$, therefore we only need to consider the cases $r\ge 0$. For $r\ge 0$, we see that $a_{|n-r|}, b_r\neq 0$ and
$$b_{r}-b_{r+\nu}>0\;\;\text{and}\;\;a_{n-r}-a_{n+r+\nu}\ge a_{|n-r|}-a_{|n-r|+ \nu}>0,$$
thus $c_n>c_{n+\nu}$. Combining the above arguments, we conclude that the sets $\mathscr{U}_x^\nu$ and $\mathscr{T}_x^\nu$ are both closed under the multiplication of Laurent polynomials. This completes the proof.
\end{proof}
Recall that a \emph{formal bivariate Laurent series} in the variables $x$ and $y$ over $\br$ has a formal expansion of the form
$\sum_{i,j\in\bz}c_{ij}x^iy^j$,
where at most finitely many $c_{ij}$ with $i<0$ or $j<0$ are non-zero.
Then, Proposition \ref{main0} can be extended into a more powerful version of formal bivariate Laurent series. To do this, we introduce the following definition.
\begin{definition}\label{defm0}
For each $\nu\in\{1,2\}$, define $\mathscr{U}_{x,y}^\nu$ and $\mathscr{T}_{x,y}^\nu$ as the sets of all formal bivariate Laurent series $f(x,y)=\sum_{i,j\in\bz}c_{ij}x^iy^j$ such that
$\sum_{i\in\bz}c_{ij}x^i\in \mathscr{U}_{x}^\nu$ and $\sum_{i\in\bz}c_{ij}x^i\in \mathscr{T}_{x}^\nu$, respectively, for any $j\in\bz$.
\end{definition}
We make the following remarks.
\begin{remark}\label{remark}By Definitions \ref{defm} and \ref{defm0}:
\begin{enumerate}
\item It is clear that $\mathscr{U}_{x,y}^1\subseteq \mathscr{U}_{x,y}^2$, $\mathscr{T}_{x,y}^1\subseteq \mathscr{T}_{x,y}^2$ and $\mathscr{T}_{x,y}^\nu\subset \mathscr{U}_{x,y}^\nu$ for each $\nu\in\{1,2\}$.
\item Any formal Laurent series $f(y)=\sum_{j}c_jy^j$ with nonnegative coefficients is an element of all $\mathscr{U}_{x,y}^\nu$ and $\mathscr{T}_{x,y}^\nu$, where $\nu\in\{1,2\}$.
\end{enumerate}
\end{remark}

As a direct consequence, we have the following extended version of Proposition \ref{main0}, which is the heart of this paper.
\begin{thm}\label{main1}For each $\nu\in\{1,2\}$, the sets $\mathscr{U}_{x,y}^\nu$ and $\mathscr{T}_{x,y}^\nu$ both form commutative semirings under the operations of addition and multiplication of formal Laurent series.
\end{thm}
\begin{proof}
By definition, it is clear that the sets $\mathscr{U}_{x,y}^\nu$ and $\mathscr{T}_{x,y}^\nu$ are both closed under the addition of formal Laurent series. Moreover, for $$f(x,y)=\sum_{i,j\in\bz}a_{ij}x^iy^j\;\;\text{and}\;\; g(x,y)=\sum_{i,j\in\bz}b_{ij}x^iy^j,$$
we have
$$f(x,y)g(x,y)=\sum_{j\in\bz}y^j\sum_{j_1+j_2=j}\bigg(\sum_{i,j_1\in\bz}a_{ij_1}x^i\bigg)\bigg(\sum_{i,j_1\in\bz}b_{ij_2}x^i\bigg)=:\sum_{j\in\bz}C_j(x)y^j.$$
By Definition \ref{defm0} and Proposition \ref{main0}, we have if $f,g\in \mathscr{U}_{x,y}^\nu$ then $C_j(x)\in \mathscr{U}_x^\nu$ for any $j\in\bz$; and if $f,g\in \mathscr{T}_{x,y}^\nu$ then $C_j(x)\in \mathscr{T}_x^\nu$ for any $j\in\bz$. Thus by Definition \ref{defm0}, we immediately obtain the proof of the theorem.
\end{proof}
In the following four sections, we will present some applications of Theorem \ref{main1}, which provides the initial motivation for this paper.

In Section~\ref{sec2}, we settle an open problem posed by Andrews \cite{MR743546} regarding the unimodality of generalized Gauss polynomials, and then we establish the unimodality of refined plane partitions studied by Okounkov and Reshetikhin \cite{MR2276355}.

In Section~\ref{sec3}, we first provide a criterion for strict unimodality induced by unimodality. Then, we turn to the unimodality and strict unimodality of certain refined color partitions.

In Section~\ref{sec4}, we study the unimodality and strict unimodality of the ranks of several combinatorial objects that generalize usual partitions, particularly focusing on concave and convex compositions studied by Andrews \cite{MR743546}, and certain unimodal sequences studied by Kim-Lim-Lovejoy \cite{MR3513531}.

Finally, we close in Section~\ref{sec5} by examining the unimodality and strict unimodality of the Betti numbers and Gromov-Witten invariants of certain Hilbert schemes of points related to integer partitions, as discussed in G\"{o}ttsche \cite{MR1032930} and Oberdieck \cite{MR3720346}.
\subsubsection*{Notations}
Let $\bn$ and $\bn_0$ be the sets of positive and non-negative
integers, respectively. Moreover, we denote that $\bn_{\infty}=\bn_0\cup\{\infty\}$.
Throughout this paper, the \textit{$q$-shifted factorial}
is defined by
$$
(a;q)_n:=\prod_{0\le j<n}(1-aq^j)
$$
for any indeterminant $a$ and $n\in \bn_{\infty}$. Moreover, the products of $q$-shifted factorials are compactly denoted as
$$
(a_1,a_2,\ldots, a_m;q)_n:=\prod_{1\le j\le m}(a_j;q)_n
$$
for any integer $m\ge 1$, indeterminant $a_1,a_2,\ldots,a_m$ and $n\in \bn_{\infty}$.

\section{Unimodality and Gauss polynomials}\label{sec2}

One of famous example of the symmetric unimodal polynomials is the
Gaussian polynomial. Let $m,n$ be any non-negative integers, then the
Gaussian polynomial $\qbm{n}{m}_q$ is defined by
\begin{align*}
\qbm{n}{m}_q=\begin{cases}
\frac{(q;q)_n}{(q;q)_m (q;q)_{n-m}}\qquad &\text{if}~0\le m\le n,\\
\qquad 0 &\text{otherwise}.
\end{cases}
\end{align*}

We begin with the following lemma, which will be used in the following discussion. It extends and rewrites the unimodality of the Gaussian polynomials into a bivariate Laurent series form as given in Definition \ref{defm0}.

\begin{lemma}\label{lemm10}For any $n\in\bn_0$, we have
\begin{align*}
(-zq^{1-n};q^2)_n\in \mathscr{U}_{q,z}^2
\quad
\text{and}
\quad
\frac{1}{(zq^{1-n};q^2)_n}\in \mathscr{U}_{q,z}^2.
\end{align*}
\end{lemma}
\begin{proof}
The use of Cauchy's binomial theorem (cf.\ Andrews \cite[Theorem 3.3]{MR1634067}) implies that
\begin{align*}
(-zq^{1-n};q^2)_n&=\sum_{0\le \ell\le n}\qbm{n}{\ell}_{q^2}(zq^{-n-1})^\ell q^{\ell(\ell+1)}=\sum_{0\le \ell\le n}q^{-(n-\ell)\ell}\qbm{n}{\ell}_{q^2}z^\ell,
\end{align*}
and
\begin{align*}
\frac{1}{(zq^{1-n};q^2)_n}&=\sum_{\ell\ge 0}\qbm{n-1+\ell}{\ell}_{q^2}(zq^{1-n})^\ell=\sum_{\ell\ge 0}q^{-(n-1)\ell}\qbm{n-1+\ell}{\ell}_{q^2}z^\ell.
\end{align*}
Notice that $q^{-(n-\ell)\ell}\qbm{n}{\ell}_{q^2}, q^{-(n-1)\ell}\qbm{n-1+\ell}{\ell}_{q^2}\in \mathscr{U}_q^2$ holds for any $\ell, n\in\bn_0$, we obtain the proof of the lemma.
\end{proof}
In the following two subsections, we first examine the unimodality of generalized Gauss polynomials, as studied by Andrews \cite{MR743546}. In particular, we solve an open problem posed by Andrews. Then, we establish the unimodality of refined plane partitions, as studied by Okounkov and Reshetikhin \cite{MR2276355}.

\subsection{On Andrews's generalized Gauss polynomials}
Following Andrews in his 1984 AMS memoir \cite{MR743546}, we define a generalised Frobenius
partition (or $F$-partition) of $n$ to be an array of nonnegative integers of the form
$$\left(\begin{array}{cccc}
   a_1 & a_2 & \ldots & a_r \\
   b_1 & b_2 & \ldots & b_r
 \end{array}\right),
$$
where each row is arranged in nonincreasing order and $n=r+\sum_{1\le i\le r}(a_i+b_i)$. Let $p(m,n,j)$ denote the number of $F$-partitions of $j$ with distinct parts in each row wherein each top entry is $<m$ and each bottom entry is $<n$. Let ${\rm CT}_z[f(z)]$ denote constant term of the formal Laurent series $f(z)$.
Andrews \cite[p.15, Equation (6.3)]{MR743546} proved that
\begin{equation*}
\sum_{0\le j\le mn}p(m,n,j)q^n={\rm CT}_z\left[(-qz;q)_m(-z^{-1};q)_n\right]=\qbm{m+n}{m}_q,
\end{equation*}
for any $m,n\in\bn_0$. This provides an explanation of the generalized Frobenius partition terminology for the Gauss polynomials. Since the Gauss polynomials have played such an important role in many branches of mathematics besides the theory of partitions, Andrews \cite[Section 6]{MR743546} introduced two types of generalized Gauss polynomials defined as
\begin{equation*}
\qbml{}{}_{\rm rk}{\qbm{m+n}{m}}_q={\rm CT}_z\left[\frac{(q^{k+1}z^{k+1};q^{k+1})_m(z^{-k-1};q^{k+1})_n}{(qz;q)_m(z^{-1};q)_n}\right]
\end{equation*}
and
\begin{equation*}
\qbml{}{}_{\rm ck}{\qbm{m+n}{m}}_q={\rm CT}_z\left[(-qz;q)_m^k(-z^{-1};q)_n^k\right],
\end{equation*}
for any integers $k\ge 1$ and $m,n\ge 0$. It is clear that for $k=1$,
$$\qbml{}{}_{\rm r1}{\qbm{m+n}{m}}_q=\qbml{}{}_{\rm c1}{\qbm{m+n}{m}}_q={\qbm{m+n}{m}}_q,$$
is the original Gauss polynomials.

\medskip

In the concluding section of the Memoir, Andrews \cite[Section 15, Problem 3]{MR743546} presented the following problem.
\begin{problem}Show that the generalized Gauss polynomials
$
\qbml{}{}_{\rm rk}{\qbm{m+n}{m}}_q
$ and $
\qbml{}{}_{\rm ck}{\qbm{m+n}{m}}_q
$
are unimodal for all integers $m,n\ge 0$ and $k\ge 1$.
\end{problem}
To solve the above problem, we define for any integers $n\ge 0$ and $d\ge 1$ that
\begin{align}
G_{d,n}(q,z):=\frac{(z^{d+1}q^{(d+1)(1-n)};q^{2(d+1)})_n}{(zq^{1-n};q^2)_n}.
\end{align}
Then, Lemma \ref{lemm10} stated that $G_{1,n}(q,z)\in\mathscr{U}_{q,z}^2$ and $$G_{\infty,n}(q,z):=\frac{1}{(zq^{1-n};q^2)_n}\in \mathscr{U}_{q,z}^2.$$ On the other hand, using the definition of the $q$-shifted factorial gives
\begin{align*}
(qz;q)_m(z^{-1};q)_n&=(qz;q)_m\prod_{0\le j<n}(1-z^{-1}q^j)\\
&=(qz;q)_m (-z)^{-n}q^{n(n-1)/2}\prod_{0\le j< n}(1-zq^{1-n+j})\\
&=(-z)^{-n}q^{n(n-1)/2}(zq^{1-n};q)_n(qz;q)_m\\
&=(-z)^{-n}q^{n(n-1)/2}(zq^{1-n};q)_{m+n},
\end{align*}
that is
$$(q^2z;q^2)_m(z^{-1};q^2)_n=(-z)^{-n}q^{n(n-1)}((zq^{1+m-n})q^{1-m-n};q^2)_{m+n}.$$
Therefore,
\begin{equation}\label{eqggpv}
\frac{(q^{2(d+1)}z^{d+1};q^{2(d+1)})_m(z^{-d-1};q^{2(d+1)})_n}{(q^2z;q^2)_m(z^{-1};q^2)_n}
=z^{-dn}q^{dn(n-1)}G_{d,m+n}(q,zq^{1+m-n}).
\end{equation}
Letting $d=1$ in the above and using Lemma \ref{lemm10} and Theorem \ref{main1}, we immediately establish the unimodality for $\qbml{}{}_{\rm ck}{\qbm{m+n}{m}}_{q}$, because of
\begin{equation*}
\qbml{}{}_{\rm ck}{\qbm{m+n}{m}}_{q^2}={\rm CT}_z\left[(-q^2z;q^2)_m^k(-z^{-1};q^2)_n^k\right]={\rm CT}_z\left[z^{-kn}q^{kn(n-1)}G_{1,m+n}(q,zq^{1+m-n})^k\right].
\end{equation*}
\begin{theorem}
For all integers $m,n\ge 0$ and $k\ge 1$, the Andrews generalized Gauss polynomials
$
\qbml{}{}_{\rm ck}{\qbm{m+n}{m}}_q
$
are unimodal.
\end{theorem}
Moreover, by \eqref{eqggpv}, we see that
\begin{equation*}
\qbml{}{}_{\rm rk}{\qbm{m+n}{m}}_{q^2}={\rm CT}_z\left[z^{-kn}q^{kn(n-1)}G_{k,m+n}(q,zq^{1+m-n})\right].
\end{equation*}
Therefore, if $G_{k,n}(q,z)\in\mathscr{U}_{q,z}^2$ for all $n\ge 0$ then the unimodality for $\qbml{}{}_{\rm rk}{\qbm{m+n}{m}}_{q}$ will follow. However, the numerical data only support the following conjecture. Using {\bf Mathematica}, we see that this conjecture is true for $k\in\{3,5,7,9,11\}$ and $n\le 10$.
\begin{conjecture}
For any odd integer $k>1$, we have $G_{k,n}(q,z)\in \mathscr{U}_{q,z}^2$ holds for all $n\ge 0$.
\end{conjecture}
Our current approach doesn't work to solve the unimodality for $\qbml{}{}_{\rm rk}{\qbm{m+n}{m}}_{q}$.
We only have a small result on this direction. Note that
$$G_{3,n}(q,z)=\frac{(z^{4}q^{4(1-n)};q^{8})_{n}}{(z^2q^{2(1-n)};q^4)_{n}}\frac{(z^{2}q^{2(1-n)};q^{4})_{n}}{(zq^{1-n};q^2)_{n}}=G_{1,n}(q^2,z^2)G_{1,n}(q,z),$$
and $G_{1,2n+1}(q,z)\in \mathscr{U}_{q^2,z}^1\subseteq \mathscr{U}_{q^2,z}^2$. We have
$$G_{3,2n+1}(q,z)=G_{1,2n+1}(q^2,z^2)G_{1,2n+1}(q,z)\in\mathscr{U}_{q^2,z}^2.$$
This immediately yields the following \emph{parity unimodality}.
\begin{proposition}For all positive integers $m\neq n\pmod 2$, we have
$$
\qbml{}{}_{\rm r3}{\qbm{m+n}{m}}_{q}q^{-\frac{3mn}{2}}\in \mathscr{U}_{q}^2.
$$
\end{proposition}
\subsection{On a refined plane partition} A plane partition $\pi$ is given by a two dimensional array of nonnegative integers $\pi_{i,j}$ that are weakly decreasing in both indices. This means that $\pi_{i,j}\ge \pi_{i+a, j+b}$ for all $i,j\in\bn$ and $a,b\in\bn_0$. Let $\mathrm{PP}$ denote the set of all plane partitions. For a plane partition $\pi=(\pi_{i,j})$, define the statistics
$$w_-(\pi)=\sum_{i<j}\pi_{i,j},\;\; w_0(\pi)=\sum_{i=j}\pi_{i,j},\;\; w_+(\pi)=\sum_{i>j}\pi_{i,j},$$
and $w(\pi)=w_-(\pi)+w_0(\pi)+w_+(\pi)$, the total sum of the integers $\pi_{i,j}$ of the partition $\pi$. For
a half-integer $\delta\in \frac{1}{2}\bz$, Okounkov and Reshetikhin established the generating function \cite{MR2276355} (see also \cite[p.151, Theorem A.1]{MR3032328}) as follows.
\begin{align}\label{eqppg}
\sum_{\pi\in \mathrm{PP}}t_1^{w_-(\pi)+(\frac{1}{2}+\delta)w_0(\pi)}t_2^{w_+(\pi)+(\frac{1}{2}-\delta)w_0(\pi)}=\prod_{i,j\ge 1}\left(1-t_1^{i-\frac{1}{2}+\delta}t_2^{j-\frac{1}{2}-\delta}\right)^{-1}.
\end{align}
Let $t_1=t_2=t$ in \eqref{eqppg}, then we obtain the well--known generating function for plane partitions:
\begin{align*}
\sum_{\pi\in \mathrm{PP}}t^{w(\pi)}=\prod_{i,j\ge 1}\left(1-t^{i+j-1}\right)^{-1}=\prod_{m\ge 1}\frac{1}{(1-t^m)^m}.
\end{align*}
For any integers $n\ge 0$ and $\delta$, we define the \emph{refined plane partition function} $b_\delta(\ell,n)$ as
$$b_\delta(\ell, n)=\#\{\pi\in \mathrm{PP}: w(\pi)=n, w_-(\pi)-w_+(\pi)+\delta w_0(\pi)=\ell\}.$$
\begin{remark}Following Morrison \cite{MR3248057}, we see that $b_\delta(\ell, n)$ are the Betti numbers of the cohomological Hall algebra, that is the refined Donaldson-Thomas invariants of $\mathbb{C}^3$. Morrison \cite{MR3248057} proved that under suitably normalized, the Betti numbers $b_\delta(\ell,n)$ have a Gaussian distribution as limit law, i.e., for large $n$ plotting the Betti numbers $b_\delta(\ell,n)$ against cohomological degree $\ell$ should give the bell curve of a Gaussian distribution.
\end{remark}
Setting $t_1=qt$ and $t_2=q^{-1}t$ in \eqref{eqppg}, we obtain
\begin{equation*}
\sum_{n\ge 0}\sum_{\ell\in \bz}b_\delta(\ell, n)q^\ell t^n=\prod_{i,j\ge 1}\frac{1}{1-q^{j-i+\delta}t^{i+j-1}}=\prod_{m\ge 1}\prod_{0\le k<m}\frac{1}{1-q^{2k+1-m+\delta}t^m}.
\end{equation*}
By Lemma \ref{lemm10}, we see that $(t^mq^{1-m};q^2)_m^{-1}\in \mathscr{U}_{q, t}^2$ for any integer $m\ge 1$. Therefore, the use of Theorem \ref{main1} immediately implies
$$\sum_{n\ge 0}\sum_{\ell\in \bz}b_0(\ell, n)q^\ell t^n=\prod_{m\ge 1}\frac{1}{(t^mq^{1-m};q^2)_m}\in \mathscr{U}_{q, t}^2.$$
Thus we immediately obtain the following parity unimodality for $b_0(m,n)$.
\begin{theorem}For all integers $m,n\ge 0$, we have
$b_0(m,n)\ge b_0(m+2,n)$.
\end{theorem}
Based on the numerical data, we propose the following conjecture. Using {\bf Mathematica}, we see that this conjecture is true for all integers $n\le 30$.
\begin{conjecture}
The sequences $(b_{0}(m, n))_{m \in \mathbb{Z}}$ are unimodal for all integers $n \geq 0$.
\end{conjecture}

\section{Unimodality and certain refined color partitions}\label{sec3}
In this section, we examine the unimodality and strict unimodality of certain \emph{refined color partitions}.
To achieve this, we first examine a class of formal bivariate Laurent series $f(z, q) \in \mathscr{U}_{z,q}^\nu$ such that $(1-q)^{-1}f(z,q)\in \mathscr{T}_{z, q}^\nu$.
\subsection{Unimodality and strict unimodality}
Let integers $\nu\in\{1,2\}$ and $0\le t<\nu$. By definition, for any
\begin{align}\label{defsssf}
f(z,q)=\sum_{n\ge 0}\sum_{i\in\bz}c_{i,n}z^iq^n\in \mathscr{U}_{z,q}^\nu,
\end{align}
and all integers $i\equiv t\pmod \nu$, there exists an integer-valued sequence $\left(\ell_n(\nu, t)\right)_{n\ge 0}$ such that the coefficients $c_{i,n}$ are positive for all $|i|\le \ell_n(\nu, t) $, and  $c_{i,n}$ equals $0$ for all $|i|>\ell_n(\nu,t)$.
The following proposition provides a criterion for strict unimodality induced by unimodality.
\begin{proposition}\label{prop321}Let $f(z,q)$ be defined by \eqref{defsssf}, let
$S_{i,n}(\nu,t)=\{0\le r\le n: \ell_r(\nu, t)\ge |i|\}$ and let $m_n(\nu,t)=\max\limits_{0\le r\le n}\ell_r(\nu, t)$. Then, we have $(1-q)^{-1}f(z,q)\in \mathscr{T}_{z,q}^\nu$, provided that
$$
\max_{r\in S_{i,n}(\nu,t)}\left(c_{i,r}-c_{i+\nu,r}\right)>0
,$$
holds for all integers $0\le t<\nu$, and $i,n\ge 0$ such that $i\le m_n(\nu,t)$ with $\ell_n(\nu,t)\ge 0$.

\end{proposition}
\begin{proof}Note that
\begin{align*}
(1-q)^{-1}f(z,q)=\sum_{j\ge 0}q^j\sum_{r\ge 0}q^r\sum_{i\in\bz}c_{i,r}z^i=\sum_{n\ge 0}q^n\sum_{i\in\bz}z^i\sum_{0\le r\le n}c_{i,r}.
\end{align*}
Suppose now that $i\equiv t\pmod 2$. Then,
$$\sum_{0\le r\le n}c_{i,r}=\sum_{r\in S_{i,n}(\nu,t)}c_{i,r}.$$
Note that if $|i|>m_n(\nu,t)$ then $S_{i,n}(\nu,t)=\varnothing$. Thus $\sum_{0\le r\le n}c_{i,r}=0$ for all $|i|>m_n(\nu,t)$. By definition, if $|i|\le m_n(\nu,t)$ then $S_{i,n}(\nu,t)\neq \varnothing$.
Therefore, by the definition of $\mathscr{T}_{z,q}^\nu$, it is sufficient to show that
$$\Delta_{i,n}^{\nu,t}:=\sum_{0\le r\le n}c_{i,r}-\sum_{0\le r\le n}c_{i+\nu,r}=\sum_{r\in S_{i,n}(\nu,t)}\left(c_{i,r}-c_{i+\nu,r}\right)>0,$$
for all integers $0\le t<\nu$, and $i,n\ge 0$ such that $i\le m_n(\nu,t)$ with $\ell_n(\nu,t)\ge 0$. In this case, the fact that $f(z,q)\in\mathscr{U}_{z,q}^\nu$ implies $c_{i,n}-c_{i+\nu,n}\ge 0$ for all $i\ge 0$. Therefore
$$
\Delta_{i,n}^\nu\ge \max_{r\in S_{i,n}(\nu,t)}\left(c_{i,r}-c_{i+\nu,r}\right)>0,$$
by using the condition of the proposition. This completes the proof of the proposition.
\end{proof}

The following corollary is a special case of Proposition \ref{prop321}, which is convenient and useful in this paper.
\begin{corollary}\label{prop3210}For $\ell_n(\nu,t)=\min(\gamma n, 2\delta-\gamma n)$ or $\ell_n(\nu,t)=\min(\gamma n, \delta)$, where $\gamma:=\gamma_{\nu,t}\in \bn$ and $\delta:=\delta_{\nu,t}\in \frac{1}{2}\bn_{\infty}$, if
$c_{i,\lceil i/\gamma\rceil}>c_{i+\nu,\lceil i/\gamma\rceil}$ holds for all $0\le i\le \delta$ and $0\le t<\nu$, then
we have
\begin{equation}\label{eqfbsu}
(1-q)^{-1}f(z,q)\in \mathscr{T}_{z,q}^\nu.
\end{equation}
In particular, if $\gamma\in\{1,\nu\}$ then \eqref{eqfbsu} holds.
\end{corollary}
\begin{proof}
For the case $\ell_n(\nu,t)=\min(\gamma n, 2\delta-\gamma n)$, note that
$$S_{i,n}(\nu,t)=\{0\le r\le n: \min(\gamma r, 2\delta-\gamma r)\ge i\}=\{r\in\bn_0: i/\gamma\le r\le \min(n, (2\delta-i)/\gamma )\}\neq \varnothing,$$
for all $i\equiv t\pmod \nu$, $0\le i\le \max\limits_{0\le r\le n}\min(\gamma r, 2\delta-\gamma r)\le \delta$ with $0\le n\le 2\delta/\gamma$. Therefore, by the condition we find that
\begin{align}\label{eqscc}
\max_{r\in S_{i,n}(\nu,t)}\left(c_{i,r}-c_{i+\nu,r}\right)\ge c_{i,\lceil i/\gamma\rceil}-c_{i+\nu,\lceil i/\gamma\rceil}>0,
\end{align}
holds for all integers $0\le t<\nu$, $0\le i\le \delta$ and $0\le n\le 2\delta/\gamma$.
Thus by Proposition \ref{prop321}, we find that $(1-q)^{-1}f(z,q)\in \mathscr{T}_{z,q}^\nu$. For $\gamma\in\{1,\nu\}$, note that
$$i+\nu>\gamma\lceil i/\gamma\rceil\ge \min(\gamma \lceil i/\gamma\rceil, 2\delta-\gamma \lceil i/\gamma\rceil)=\ell_{\lceil i/\gamma\rceil}(\nu,t),$$
we see that $c_{i+\nu,\lceil i/\gamma\rceil}=0$, and thus \eqref{eqscc} is true. Therefore if $\gamma\in\{1,\nu\}$ then \eqref{eqfbsu} holds.

The proof for the case $\ell_n(\nu,t)=\min(\gamma n, \delta)$ is similar, we omit the details. This completes the proof of the corollary.
\end{proof}

\subsection{On certain refined color partitions}
Recall that a partition $\lambda$ is a non-increasing finite sequence $\lambda=(\lambda_1,\lambda_2,\ldots,\lambda_s)$ of positive integers. The elements $\lambda_i$ that appear in the sequence $\lambda$ are called parts of $\lambda$. We say $\lambda$ is a partition
of $n$, if the sum $|\lambda|:=\sum_{i}\lambda_i$ of all parts of $\lambda$ is equal to $n$ (see the reference by Andrews \cite{MR1634067}).

\medskip

Let $S$ be \emph{any set (could be a multiset) of positive integers}. For any $a\in\{1,2\}$ and $b\in \bn_\infty$, we define the (formal) Laurent series for $\mathfrak{S}_{ab}(z,q)$ by
\begin{align}\label{znhp}
\mathfrak{S}_{ab}(z,q)&:=\sum_{n\ge 0}\sum_{m\in\bz}p_{ab,S}(m,n)z^mq^n\nonumber\\
&=\prod_{\ell\in S}\frac{(1-q^{(3-a)\ell})\left(1-(zq^\ell)^{b+1}\right)\left(1-(z^{-1}q^\ell)^{b+1}\right)}{(1-q^\ell)(1-zq^\ell)(1-z^{-1}q^\ell)}.
\end{align}
Moreover, we define that
\begin{align*}
\mathfrak{S}_{ab}(q):=\sum_{n\ge 0}p_{ab,S}(n)q^n=\prod_{\ell\in S}\frac{(1-q^{(3-a)\ell})\left(1-q^{(b+1)\ell}\right)^2}{(1-q^\ell)(1-q^\ell)^2}.
\end{align*}
Then, the coefficient $p_{ab,S}(n)$ in the power series for $\mathfrak{S}_{ab}(q)$ counts the number of partitions of $n$ into parts of three colors, say, red, green, and blue, such that all of the red, green, and blue parts are chosen from $S$, subject to the restriction that each of the red parts appears at most $(2-a)$ times, and each of green and blue parts appears at most $b$ times. Furthermore, the coefficient $p_{ab,S}(m,n)$ in the Laurent series for $\mathfrak{S}_{ab}(z, q)$ counts the number of partitions of $n$ that are same to those for $p_{ab,S}(n)$, with the additional restriction that \emph{the difference between the numbers of green and blue parts equals $m$}. Therefore, $p_{ab,S}(m,n)$ refines $p_{ab,S}(n)$.

\medskip

Similarly to $\mathfrak{S}_{ab}(z,q)$, we define the (formal) Laurent series for $\mathbf{S}_{ab}(z,q)$ with $a\in\{1,2\}$ and $b\in \bn_\infty$ by
\begin{align}\label{znhb}
\mathbf{S}_{ab}(z,q)&:=\sum_{n\ge 0}\sum_{m\in\bz}r_{ab,S}(m,n)z^mq^n\nonumber\\
&=\prod_{\ell\in S}\frac{(1-q^{(3-a)\ell})\left(1-(zq^\ell)^{b+1}\right)\left(1-(z^{-1}q^\ell)^{b+1}\right)}{(1-q^\ell)^2(1-zq^\ell)(1-z^{-1}q^\ell)}.
\end{align}
Moreover, we define that
\begin{align*}
\mathbf{S}_{ab}(q):=\sum_{n\ge 0}r_{ab,S}(n)q^n=\prod_{\ell\in S}\frac{(1-q^{(3-a)\ell})\left(1-q^{(b+1)\ell}\right)^2}{(1-q^\ell)(1-q^\ell)(1-q^\ell)^2}.
\end{align*}
It is clear that the Laurent series for $\mathbf{S}_{ab}(z, q)$ defined by \eqref{znhb} also generate certain refined color partitions, similar to $\mathfrak{S}_{ab}(z,q)$. In particular, the coefficient $r_{ab,S}(n)$ in the power series for $\mathbf{S}_{ab}(q)$ counts the number of partitions of $n$ into parts of four colors, say, yellow, red, green, and blue, such that all of the yellow, red, green, and blue parts are chosen from $S$, subject to the restriction that each of the red parts appears at most $(2-a)$ times, and each of green and blue parts appears at most $b$ times. Meanwhile, the coefficient $r_{ab,S}(m,n)$ in the Laurent series for $\mathbf{S}_{ab}(z, q)$ counts the number of partitions of $n$ that are same to those for $r_{ab,S}(n)$, with the additional restriction that \emph{the difference between the numbers of green and blue parts equals $m$}. That is $r_{ab,S}(m,n)$ refines $r_{ab,S}(n)$.

\medskip

Based on Theorem \ref{main1}, we prove the following theorem.
\begin{theorem}\label{main}
For any $a\in\{1,2\}$ and $b\in \bn_{\infty}$, we have $\mathfrak{S}_{ab}(z,q)\in \mathscr{U}_{z,q}^a$ and $\mathbf{S}_{ab}(z,q)\in \mathscr{T}_{z,q}^a$.
\end{theorem}
\begin{remark}
Let $S=\bn$, $a\in\{1,2\}$ and $b\in \bn_{\infty}$. Then we have the following special cases:
\begin{equation*}
\frac{(q^{3-a};q^{3-a})_\infty\left(z^{b+1}q^{b+1},z^{-b-1}q^{b+1};q^{b+1}\right)_\infty}{(q;q)_\infty(zq,z^{-1}q;q)_\infty}\in\mathscr{U}_{z,q}^a,
\end{equation*}
and
\begin{equation*}
\frac{(q^{3-a};q^{3-a})_\infty\left(z^{b+1}q^{b+1},z^{-b-1}q^{b+1};q^{b+1}\right)_\infty}{(q;q)_\infty^2(zq,z^{-1}q;q)_\infty}\in\mathscr{T}_{z,q}^a.
\end{equation*}
\end{remark}

As a direct consequence of Theorem \ref{main} and Definition \ref{defm0}, we have the following monotonicity results for $p_{ab,S}(m,n)$ and $r_{ab,S}(m,n)$.
\begin{theorem}
Defined $p_{ab, S}(m,n)$ and $r_{ab,S}(m,n)$ by \eqref{znhp} and \eqref{znhb}, respectively.
For any $b\in\bn_\infty$, any integers $a\in\{1,2\}$ and $m,n\ge 0$, we have $p_{ab, S}(m,n)\ge p_{ab, S}(m+a,n)$; and
\begin{align*}
r_{ab, S}(m,n)>r_{ab, S}(m+a,n),\;\text{provided}\; r_{ab,S}(m,n)\neq 0.
\end{align*}
\end{theorem}
\begin{remark}
By \eqref{znhp} and \eqref{znhb}, we see that if $S$ only contains positive odd integers, then for all integers $m \not\equiv n \pmod{2}$, we have $p_{2b,S}(m,n)=r_{2b,S}(m,n)=0$.
\end{remark}

It is not difficult to see that the proof of Theorem \ref{main} follows from Theorem \ref{main1} and the following Lemma \ref{lemm11}.
\begin{lemma}\label{lemm11}Let $a\in\{1,2\}$ and $b\in\bn_\infty$. We have
\begin{equation}\label{eqim1}
\frac{(1-q^{3-a})\left(1-(zq)^{b+1}\right)\left(1-(z^{-1}q)^{b+1}\right)}{(1-q)(1-zq)(1-z^{-1}q)}\in \mathscr{U}_{z,q}^a,
\end{equation}
and
\begin{equation}\label{eqim10}
\frac{(1-q^{3-a})\left(1-(zq)^{b+1}\right)\left(1-(z^{-1}q)^{b+1}\right)}{(1-q)^2(1-zq)(1-z^{-1}q)}\in \mathscr{T}_{z,q}^a.
\end{equation}
\end{lemma}
\begin{proof}
Note that
\begin{align*}
\frac{\left(1-(zq)^{b+1}\right)\left(1-(z^{-1}q)^{b+1}\right)}{(1-zq)(1-z^{-1}q)}=\sum_{0\le r,s\le b}q^{r+s}z^{r-s}=\sum_{k=0}^{2b}q^{k}\sum_{\substack{r+s=k\\ 0\le r,s\le b}}z^{r-s}=\sum_{0\le k\le 2b}q^k\sum_{\substack{|r|\le \min(k,2b-k)\\ r\equiv k\pmod 2}}z^{r}.
\end{align*}
We see that the relation \eqref{eqim1} is true for $a=2$. Furthermore, by using Corollary \ref{prop3210}, we find that \eqref{eqim10} is true for $a=2$. For $a=1$, note that
\begin{align*}
\frac{(1+q)\left(1-(zq)^{b+1}\right)\left(1-(z^{-1}q)^{b+1}\right)}{(1-zq)(1-z^{-1}q)}=&(1+q)\sum_{0\le k\le 2b}q^k\sum_{\substack{|r|\le \min(k,2b-k)\\ r\equiv k\pmod 2}}z^{r}\\
=&1+q^{2b+1}+\sum_{1\le k\le b}q^k\bigg(\sum_{\substack{|r|\le \min(k,2b-k)\\ r\equiv k\pmod 2}}z^{r}+\sum_{\substack{|r|\le \min(k-1,2b+1-k)\\ r\equiv k-1\pmod 2}}z^{r}\bigg).
\end{align*}
That is
\begin{align}\label{eqmmm0}
\frac{(1+q)\left(1-(xq)^{b+1}\right)\left(1-(x^{-1}q)^{b+1}\right)}{(1-xq)(1-x^{-1}q)}=\sum_{k=0}^{2b+1}q^k\sum_{|r|\le \min(k, 2b+1-k)}x^{r}.
\end{align}
Therefore, the relation \eqref{eqim1} is true for $a=1$. Furthermore, by using Corollary \ref{prop3210}, we find that \eqref{eqim10} is true for $a=1$, which completes the proof.
\end{proof}
As the end of this section, we define the (formal) Laurent series for $\mathsf{S}(z,q)$ by
\begin{align*}\label{znhd}
\mathsf{S}(z,q):=\sum_{n\ge 0}\sum_{m\in\bz}d_{S}(m,n)z^mq^n=\prod_{\ell\in S}\left(1+zq^\ell\right)\left(1+z^{-1}q^\ell\right).
\end{align*}
In view of
$$(1+zq)(1+z^{-1}q)=1+q^2+q(z+z^{-1})\in \mathscr{T}_{z,q}^2,$$
using Theorem \ref{main1} implies the following theorem.
\begin{theorem}\label{main00}
We have $\mathsf{S}(z,q)\in \mathscr{T}_{z,q}^2$. In particular, $d_{S}(m,n)>d_{S}(m+2,n)$, provided that $m,n\ge 0$ such that
$d_S(m,n)\neq 0$.
\end{theorem}
\section{Unimodality and ranks in partition theory}\label{sec4}
Let $p(n)$ denote the number of partitions of an integer $n$. To explain the famous partition
congruences with modulus $5$, $7$ and $11$ of Ramanujan \cite{MR2280868}, the rank and crank statistic for integer partitions was introduced and investigated by Dyson \cite{MR3077150}, Andrews and Garvan \cite{MR929094, MR920146}. As a precise definition of rank and crank for integer partitions are not necessary for the rest of the paper, we do not give it here. Let $M(m, n)$ (with a slight modification in the case that $n = 1$, where the values are instead $M(\pm 1, 1) = 1, M(0, 1) = -1$) and $N(m,n)$ denote the number of partitions of $n$ with crank $m$ and rank $m$, respectively. It is well--known that
\begin{equation*}
\mathcal{C}(x,q):=\sum_{n\ge 0}\sum_{m\in\bz}M(m,n)x^mq^n=\frac{(q;q)_\infty}{(xq,x^{-1}q;q)_\infty},
\end{equation*}
and
\begin{equation*}
\mathcal{R}(x,q):=\sum_{n\ge 0}\sum_{m\in\bz}N(m,n)x^mq^n=\sum_{n\ge 0}\frac{q^{n^2}}{(x q, x^{-1}q;q)_{n}}.
\end{equation*}

The rest of this section aims to study the unimodality and strict unimodality of the ranks of several combinatorial objects generalized from usual partitions, especially for \emph{concave and convex compositions} and certain \emph{unimodal sequences}. The exploration of ranks and cranks for usual partitions has been extended to a variety of related combinatorial objects, such as unimodal sequences, concave compositions, overpartitions, spt type partitions, among others, which play an important role in the study of integer partitions.  It is worth mentioning that the combinatorial generating functions for the above objects are closely related to Jacobi forms, mock Jacobi forms, and quantum Jacobi forms, all of which play significant roles in not only combinatorics but also modern number theory (see \cite{MR2308850, MR3048655, MR3152010, MR2177037, MR2630043,MR2994899, MR3513531, MR2176595,MR2401139} for example).

\subsection{Garvan $k$-ranks}In \cite{MR1291125}, Garvan generalized Dyson's rank to the \emph{$k$-ranks} for any integer $k\ge 2$ as follows.
\begin{align*}
\mathcal{R}_k(x,q):=&\sum_{n\ge 0}\sum_{m\in\bz}N_k(m,n)x^mq^n\\
=&\sum_{n_{k-1}\ge n_{k-2}\ge \cdots\ge n_1\ge 0}\frac{q^{n_1^2+n_2^2+\cdots+n_{k-1}^2}}{(q;q)_{n_{k-1}-n_{k-2}}\cdots (q;q)_{n_3-n_2} (q;q)_{n_2-n_1}(x q, x^{-1}q;q)_{n_1}}.
\end{align*}
In particular, Garvan \cite{MR1291125} proved that for any integer $k\ge 2$, $N_k(m, n)$ is the number of partitions of $n$ into at least $(k-1)$ successive Durfee squares with \emph{$k$-rank} equal to $m$. For the detail of the combinatorial interpretation of $N_k(m,n)$, see \cite[Theorem (1.12)]{MR1291125}.
\medskip

In  \cite{MR4324846}, Ji and Zang proved that the sequences $(M(m,n))_{|m|\le n-1}$ are unimodal for all integers $n\ge 44$. However, using {\bf Mathematica}, we see that the crank generating function $\mathcal{C}(x,q)\not\in \mathscr{U}_{x,q}^2$. For the $k$-rank generating function, since (by Theorem \ref{main})
$$\frac{1}{(x q, x^{-1}q;q)_{n_1}}\in \mathscr{U}_{x,q}^2,$$
for all integers $n_1\ge 0$, the use of Theorem \ref{main1} implies that $\mathcal{R}_k(x,q)\in \mathscr{U}_{x,q}^2$ for all integers $k\ge 2$. This immediately leads to the following theorem.
\begin{theorem}For any $k\ge 2$ and $n,m\ge 0$, we have
$$N_k(m,n)\ge N_k(m+2,n).$$
\end{theorem}
Let $k=2$, then we obtain the parity unimodality of $N(m,n)=N_2(m,n)$ in $m$. This case was first proved by
Chan and Mao \cite{MR3190432} by using a combinatorial argument. The cases for $k\ge 3$ are new. We note that
Liu and Zhou \cite[Corollay 1.5]{LZ}, proved that for each positive integer $k\ge 2$, $\left(N_k(m,n)\right)_{|m|\le n-k}$ is a unimodal sequence for all sufficiently large positive integers $n$. Moreover, Liu and Zhou \cite[Theorem 1.3]{LZ} established the following uniform asymptotic formulas of $k$-rank functions:
\begin{align*}
\frac{N_k(m,n)}{p(n)}=\frac{\pi}{4\sqrt{6n}}{\rm sech}^2\left(\frac{\pi m}{2\sqrt{6n}}\right)\left(1+O\left(\frac{n+m^2}{n^{3/2}}\right)\right),
\end{align*}
as $n\to +\infty$, for each $k\ge 2$, and uniformly with respect to $m= o(n^{3/4})$.
Since $\frac{1}{4}{\rm sech}^2\left(\frac{x}{2}\right)$ is the probability density function of the standard logistic distribution, which yields that under suitably normalized, the $k$-rank function $N_k(m,n)$ have a logistic distribution as limit law.
We also note that Zhou \cite{MR4708721} has conjectured that
the sequence $\big(N_k(m,n)\big)_{|m|\le n-k}$ is a unimodal sequence for all integers $n\ge (k+1)+6\cdot {\bf 1}_{k=3}$ and $k\ge 3$.

\subsection{Ranks of concave and convex compositions}Following Andrews \cite[Section 1]{MR3048655}, a \emph{concave composition} of a positive integer $n$ is a sum of integers of the form
$$\sum_{1\le i\le R}a_i+c+\sum_{1\le i\le S}b_i=n,$$
where all of $c\ge 0$, $a_i>0$ and $b_i>0$ are called the \emph{part}, and
\begin{align}\label{eqmm1}
a_1\ge a_2\ge \dots\ge a_{R}>c<b_{1}\le \dots\le b_{S-1}\le b_S.
\end{align}
Here $c$ is called the central part. If all the ``$\le$" and ``$\ge$" in \eqref{eqmm1} are replaced by ``$<$" and ``$>$", we refer to a \emph{strictly
concave composition}. \emph{Convex and strictly convex compositions} are defined is the same way with all the
inequality signs reversed. All of the ranks of the concave composition, strictly
concave composition, convex and strictly convex compositions are defined as $S-R$, that is the number of parts after the central part minus the number of parts before the central part, which analogs the rank statistic for integer partitions and measures the position of the central part.

\medskip

Let $V(m,n), V_d(m,n), X(m,n)$ and $X_d(m,n)$ denote the numbers of concave compositions, strictly concave compositions, convex compositions and strictly convex compositions, respectively, of $n$ with rank equal to $m$. By standard counting techniques, we have the following generating functions:
\begin{align}
\label{eqv}\mathcal{V}(z,q)&=\sum_{n\ge 0}\sum_{m\in\bz}V(m,n)z^m q^n=\sum_{n\ge 0}\frac{q^{n}}{(zq^{n+1}, z^{-1}q^{n+1};q)_\infty},\\
\label{eqvd}\mathcal{V}_d(z,q)&=\sum_{n\ge 0}\sum_{m\in\bz}V_d(m,n)z^m q^n=\sum_{n\ge 0}q^{n}(-zq^{n+1}, -z^{-1}q^{n+1};q)_\infty,\\
\label{eqx}\mathcal{X}(z,q)&=\sum_{n\ge 0}\sum_{m\in\bz}X(m,n)z^m q^n=\sum_{n\ge 0}\frac{q^{n+1}}{(zq, z^{-1}q;q)_n},\\
\label{eqxd}\mathcal{X}_d(z,q)&=\sum_{n\ge 0}\sum_{m\in\bz}X_d(m,n)z^m q^n=\sum_{n\ge 0}q^{n+1}(-zq, -z^{-1}q;q)_n.
\end{align}

The generating functions for $V(m,n)$ and $V_d(m,n)$ also follow from Andrews-Rhoades-Zwegers \cite[Equations (1.7),(1.8)]{MR3152010}.\footnote{The literature \cite{MR3152010} use the word ``strongly concave composition" instead of ``strictly concave composition".} The convex compositions and strictly convex compositions are usually referred to as \emph{unimodal sequences} and \emph{strongly unimodal sequences}, respectively, in the literature; see \cite{MR2994899, MR4299082} for example.
Recall that a strongly unimodal sequence of weight $n$ is a sequence of positive integers of form
\begin{align}\label{eqmm11}
a_1<a_2<\dots< a_{R}< c> b_{1}> \dots> b_{S-1}> b_S,
\end{align}
where $\sum_{1\le i\le R}a_i+c+\sum_{1\le i\le S}b_i=n$.
The maximum value, $c$, now is called the peak.
If the inequalities are not strict, the sequence is called unimodal sequence. The rank of the above unimodal sequence is defined as $S-R$, that is the number of parts after the peak minus the number of parts before the peak. Let $u(m,n)$ and $u_w(m,n)$ denote the numbers of strongly unimodal sequences and unimodal sequences, respectively, of size $n$ with rank equal to $m$.\footnote{We note that in the related works, including \cite{MR2994899} and \cite{MR4299082}, $u_w(m,n)$ and $u(m,n)$ are both denoted as $u(m,n)$, which can lead to confusion when discussing both in the same paper.} Then, for all integers $n\ge 0$ and $m\in\bz$, from \cite[equation {[1.1]}]{MR2994899}, we have $u(m,n)=X_d(m,n)$, and from \cite[p. 455, the last line]{MR4299082}, we see that $u_w(m,n)=X(m,n+1)$.

\medskip

For the summands in \eqref{eqv}--\eqref{eqxd}, by Theorem \ref{main},
$$\frac{q^{n}}{(zq^{n+1},\; z^{-1}q^{n+1};q)_\infty}\;\;\text{and}\;\; \frac{q^{n+1}}{(zq, z^{-1}q;q)_n}$$
belong to $\mathscr{U}_{z,q}^2$, and by Theorem \ref{main00},
$$q^n(-zq^{n+1}, -z^{-1}q^{n+1};q)_\infty\;\; \text{and}\;\; q^{n+1}(-zq, -z^{-1}q;q)_n$$
belong to $\mathscr{T}_{z,q}^2$. Therefore,  the application of Theorem \ref{main1} immediately leads to the following unimodality and strict unimodality results.
\begin{theorem}\label{thumcc}We have $\mathcal{V}(z,q), \mathcal{X}(z,q) \in \mathscr{U}_{z,q}^2$.
In particular,
\begin{align*}
&V(m,n)\ge V(m+2,n)\;\;\text{and}\;\;X(m,n)\ge X(m+2,n),
\end{align*}
for all integers $m,n\ge 0$.
\end{theorem}
\begin{theorem}\label{thumcc1}We have $\mathcal{V}_d(z,q), \mathcal{X}_d(z,q)\in \mathscr{T}_{z,q}^2$.
In particular, for all integers $m,n\ge 0$,
\begin{align*}
&V_d(m,n)> V_d(m+2,n),\;\text{provided}\; V_d(m,n)\neq 0;\\
&X_d(m,n)> X_d(m+2,n),\;\text{provided}\; X_d(m,n)\neq 0.
\end{align*}
\end{theorem}
\begin{remark}By using the result from Zhou \cite[Proposition 1.2]{MR4001540}, we have:
\begin{equation*}
V_d\left(m,n+\frac{1}{2}|m|(|m|+1)\right)=\sum_{\ell\ge 0}\left(\frac{-3}{2\ell+1}\right)p\left(n-\frac{2}{3}\ell(\ell+1)-|m|\ell\right),
\end{equation*}
for any integers $m,n\in\bz$, where $(\frac{\cdot}{\cdot})$ is the Kronecker symbol.  Consequently, we have:
$$V_d(m,n)=p(n-m(m+1)/2),$$
for all integers $n\le 2m+3+m(m+1)/2$. Thus, it is not difficult to see that if $V_d(m,n)=0$ but $V_d(m,n+2)\neq 0$, then we must have $V_d(m,n+2)=1$ by using the fact that $p(0)=p(1)=1$.


\end{remark}
The numerical data suggest that Theorem \ref{thumcc1} also holds for $\mathcal{V}(z,q)$ and $\mathcal{X}(z,q)$. We have verified that Theorem \ref{thumcc1} for $V(m,n)$ and $X(m,n)$ is true for all integers $n\le 50$ using {\bf Mathematica}.
 Moreover, from the very recent work of Zang \cite[Theorem 1.3]{Zang2024v1}, we can see that $u(m,n) \ge u(m+1,n)$
for all integers $m,n \ge 0$, and $u(m,n)>u(m+1,n)$
for all $m\ge 0$ and $n\ge \max\left(6, \binom{m+2}{2}\right)$.  In particular, we can find that for all $m,n\ge 0$,
$$u(m,n)>u(m+2,n),\;  \text{provided}\; u(m,n)\neq 0.$$
In our notation, this establishes that
$\mathcal{X}_d(z,q) \in \mathscr{U}_{z,q}^1$ and $\mathcal{X}_d(z,q) \in \mathscr{T}_{z,q}^2$.
However, using {\bf Mathematica}, we observe that
$\mathcal{X}(z,q) \not\in \mathscr{U}_{z,q}^1$.
On the other hand, the numerical data suggest that a stronger conjecture regarding $\mathcal{V}(z,q)$ and $\mathcal{V}_d(z,q)$ appears to be true.
\begin{conjecture}\label{conj31}Both of $\mathcal{V}(z,q), \mathcal{V}_d(z,q)\in \mathscr{U}_{z,q}^1$. In particular,
$$V(m,n)\ge V(m+1,n)\;\;\text{and}\;\; V_d(m,n)\ge V_d(m+1,n),$$
for all integers $m,n\ge 0$.
\end{conjecture}
We have verified that Conjecture \ref{conj31} is true for all integers $n\le 50$ using {\bf Mathematica}.
It should be notice that in \cite{MR4001540}, Zhou proved that under suitably normalized, the rank function $V_d(m,n)$ of a strongly concave compositions have a Gaussian distribution as limit law, i.e., for large $n$ plotting  $V_d(m,n)$ against rank $m$ should give the bell curve of a Gaussian distribution.

\medskip

In \cite[Section 3]{MR3048655}, Andrews further explored to concave compositions and convex compositions (as defined in \eqref{eqmm1}) with the added constraint that each of the $a_i$,$b_i$ and $c$ is a positive odd integer. In these cases, we refer to them as \emph{odd concave compositions} and \emph{odd convex compositions}. Let $VO(m,n), VO_d(m,n), XO(m,n)$ and $XO_d(m,n)$ denote the numbers of odd concave compositions, odd strictly concave compositions, odd convex compositions and odd strictly convex compositions, respectively, of $n$ with rank equal to $m$. By using a similar argument for \eqref{eqv}--\eqref{eqxd}, we can derive the following generating functions:
\begin{align}
\label{eqvo}\mathcal{VO}(z,q)&=\sum_{n\ge 0}\sum_{m\in\bz}VO(m,n)z^m q^n=\sum_{n\ge 1}\frac{q^{2n-1}}{(zq^{2n+1}, z^{-1}q^{2n+1};q^2)_\infty},\\
\label{eqvod}\mathcal{VO}_d(z,q)&=\sum_{n\ge 0}\sum_{m\in\bz}{VO}_d(m,n)z^m q^n=\sum_{n\ge 1}q^{2n-1}(-zq^{2n+1}, -z^{-1}q^{2n+1};q^2)_\infty,\\
\label{eqxo}\mathcal{XO}(z,q)&=\sum_{n\ge 0}\sum_{m\in\bz}XO(m,n)z^m q^n=\sum_{n\ge 0}\frac{q^{2n+1}}{(zq, z^{-1}q;q^2)_{n}},\\
\label{eqxod}\mathcal{XO}_d(z,q)&=\sum_{n\ge 0}\sum_{m\in\bz}XO_d(m,n)z^m q^n=\sum_{n\ge 0}q^{2n+1}(-zq, -z^{-1}q;q^2)_n.
\end{align}
\begin{remark}By \eqref{eqvo}--\eqref{eqxod}, it is clear that
$$VO(m,n)=VO_d(m,n)=XO(m,n)=XO_d(m,n)=0,$$
for any integers $m,n$ with $m\not\equiv n\pmod 2$.
\end{remark}
\begin{remark}
We note that the generating functions for odd (strictly) convex compositions, as given in the above, are the essentially same as the rank generating functions for odd (strongly) unimodal sequences studied in the recent work of Bringmann and Lovejoy \cite[Theorems 1.1 and 1.2]{BL2023v1}. In particular, $XO(m,n+2)={\rm ou}(m,n)$ and $XO_d(m,n)={\rm ou}^*(m,n)$ for all integers $n\ge 0$ and $m\in\bz$, where ${\rm ou}(m,n)$ and ${\rm ou}^*(m,n)$ denote the number of odd unimodal sequences and odd strongly unimodal sequences, respectively, of weight $n$ with rank $m$.
\end{remark}
By using a similar argument for Theorems \ref{thumcc} and \ref{thumcc1}, we have the following Theorems \ref{thumocc} and \ref{thumocc1}. In particular, for the summands in \eqref{eqvo}--\eqref{eqxod}, by Theorem \ref{main},
$$\frac{q^{2n-1}}{(zq^{2n+1},\; z^{-1}q^{2n+1};q^2)_\infty}\;\;\text{and}\;\; \frac{q^{2n+1}}{(zq, z^{-1}q;q^2)_n}$$
belong to $\mathscr{U}_{z,q}^2$, and by Theorem \ref{main00},
$$q^{2n-1}(-zq^{2n+1}, -z^{-1}q^{2n+1};q^2)_\infty\;\; \text{and}\;\; q^{2n+1}(-zq, -z^{-1}q;q^2)_n$$
belong to $\mathscr{T}_{z,q}^2$. Therefore,  the application of Theorem \ref{main1} immediately leads to the following unimodality and strict unimodality results.

\begin{theorem}\label{thumocc}We have $\mathcal{VO}(z,q), \mathcal{XO}(z,q) \in \mathscr{U}_{z,q}^2$.
In particular,
\begin{align*}
VO(m,n)\ge VO(m+2,n)\;\;\text{and}\;\;XO(m,n)\ge XO(m+2,n),
\end{align*}
for all integers $m,n\ge 0$.
\end{theorem}
\begin{theorem}\label{thumocc1}We have $\mathcal{VO}_d(z,q), \mathcal{XO}_d(z,q)\in \mathscr{T}_{z,q}^2$.
In particular, for all integers $m,n\ge 0$,
\begin{align*}
&VO_d(m,n)> VO_d(m+2,n),\;\text{provided}\; VO_d(m,n)\neq 0;\\
&XO_d(m,n)> XO_d(m+2,n),\;\text{provided}\; XO_d(m,n)\neq 0.
\end{align*}
\end{theorem}
At the end of this subsubsection, we point out that the numerical data also suggest that Theorem \ref{thumocc1} also holds for $\mathcal{VO}(z,q)$ and $\mathcal{XO}(z,q)$. We have verified that Theorem \ref{thumocc1} for $VO(m,n)$ and $XO(m,n)$ is true for all integers $n\le 50$ using {\bf Mathematica}.
\subsection{Ranks of certain  unimodal sequences}
In the following we establish the unimodality of ranks of certain  unimodal sequences introduced by Kim--Lim--Lovejoy \cite{MR3513531}.
Relax the ``strong" condition in the definition (see \eqref{eqmm11}) of $u(m,n)$ and allow odd parts to repeat on either side of the peak. Then add two requirements. First the peak must be even and second, the odd parts must be identical on either side of the peak. We call these \emph{odd-balanced unimodal sequences}. Let $u_{ob}(m,n)$ denote the number of odd-balanced unimodal sequences of weight $2n+2$ and having rank $m$.\footnote{In \cite{MR3513531}, $u_{ob}(m,n)$ is denoted by $v(m,n)$ and $U_{ob}(z,q)$ is denoted by $\mathcal{V}(z,q)$.} Then by Kim--Lim--Lovejoy \cite[Equation (1.10)]{MR3513531}, we have the generating function of $u_{ob}(m, n)$ as follows
\begin{equation}\label{equsv}
U_{ob}(z,q):=\sum_{n\ge 0}\sum_{m\in\bz}u_{ob}(m,n)z^m q^n=\sum_{n\ge 0}\frac{(-zq,-z^{-1}q;q)_nq^n}{(q;q^2)_{n+1}}.
\end{equation}

Using Theorems \ref{main1}, \ref{main} and \ref{main00}, we have
\begin{align*}
&\frac{(-zq,-z^{-1}q;q)_n}{(q;q^2)_{n+1}}\in\mathscr{T}_{z,q}^2\\
=&\frac{(-zq,-z^{-1}q;q)_n(q^2;q^2)_n}{(q;q)_{2n+1}}=\frac{(-zq,-z^{-1}q,-q;q)_n(q;q)_n}{(q;q)_{2n+1}}\\
=&(-zq,-z^{-1}q,-q;q)_n\frac{1}{(q^{n+1};q)_{n+1}}\in \mathscr{U}_{z,q}^1,
\end{align*}
hence $U_{ob}(z,q)\in \mathscr{U}_{z,q}^1$ and $U_{ob}(z,q)\in\mathscr{T}_{z,q}^2$. We restated the above as the following theorem.
\begin{theorem}\label{uobus}
For all integers $n\ge 0$, the sequences $(u_{ob}(m,n))_{m\in\bz}$ defined by \eqref{equsv}, are symmetric and unimodal. Moreover, for all integers $m,n\ge 0$,
$$u_{ob}(m,n)>u_{ob}(m+2,n),$$
provided $u_{ob}(m,n)\neq 0$.
\end{theorem}

Kim--Lim--Lovejoy \cite[Section 1.3]{MR3513531} further considered two functions $\mathcal{W}(x,q)$ and $\mathcal{Z}(x,q)$, which share some, but not all, features with $\mathcal{V}(x,q)$, as follows:
\begin{equation*}
\mathcal{W}(x,q):=\sum_{n\ge 0}\frac{(xq,x^{-1}q;q^2)_nq^{2n}}{(-q;q)_{2n+1}}
\;\;
\text{and}
\;\;
\mathcal{Z}(x,q)=\sum_{n\ge 0}\frac{(-xq,-x^{-1}q;q)_nq^{n}}{(q;q)_{2n+1}}.
\end{equation*}
They \cite[Sections 3,4]{MR3513531} demonstrate that both $\mathcal{W}(x,q)$ and $\mathcal{Z}(x,q)$ have simple combinatorial interpretations as generating functions.
Define $w(m,n)$ and $z(m,n)$ as follows:
\begin{equation}\label{equsw}
\mathcal{W}(x,-q)=\sum_{n\ge 0}\sum_{m\in\mathbb{Z}}w(m,n)x^m q^n=\sum_{n\ge 0}\frac{(-xq,-x^{-1}q;q^2)_nq^{2n}}{(q;-q)_{2n+1}}
\end{equation}
and
\begin{equation}\label{equsz}
\mathcal{Z}(x,q)=\sum_{n\ge 0}\sum_{m\in\mathbb{Z}}z(m,n)x^m q^n=\sum_{n\ge 0}\frac{(-xq,-x^{-1}q;q)_nq^{n}}{(q;q)_{2n+1}}.
\end{equation}
Then $w(m,n)$ and $z(m,n)$ count certain \emph{rank statistics} of certain unimodal sequences. Notably, $z(m,n)$ counts the ranks of \emph{over-balanced unimodal sequences} (see \cite[Section 3.1]{MR3513531} for details).

\medskip

Similar to the discussion of Theorem \ref{uobus}, using Theorems \ref{main1}, \ref{main} and \ref{main00}, we have
\begin{align*}
\frac{(-xq,-x^{-1}q;q^2)_n}{(q;-q)_{2n+1}}&=\frac{(-xq,-x^{-1}q;q^2)_n}{(q;q^2)_{n+1}(-q^2;q^2)_n}\\
&=(-xq,-x^{-1}q;q^2)_n\frac{(-q, q^2;q^2)_n}{(1-q^{2n+1})(q^2;q^2)_{2n}}\in \mathscr{T}_{x,q}^2\\
&=(-xq,-x^{-1}q,-q;q^2)_n\frac{1}{(1-q^{2n+1})(q^{2n+2};q^2)_{n}}\in \mathscr{U}_{x,q}^1,
\end{align*}
and
\begin{align*}
&\frac{(-xq,-x^{-1}q;q)_n}{(q;q)_{2n+1}}\in \mathscr{T}_{x,q}^2\\
=&(-xq,-x^{-1}q,-q;q)_n\frac{1}{(q^2;q^2)_n(q^{n+1};q)_{n+1}}\in \mathscr{U}_{x,q}^1.
\end{align*}
By Theorem \ref{main1},  it follows that $\mathcal{W}(x,-q), \mathcal{Z}(x,q)\in\mathscr{U}_{x,q}^1$ and $\mathcal{W}(x,-q), \mathcal{Z}(x,q)\in\mathscr{T}_{x,q}^2$. We restated the above as the following theorem.
\begin{theorem}
For all integers $n\ge 0$, the sequences $(w(m,n))_{m\in\bz}$ and $(z(m,n))_{m\in\bz}$ defined by \eqref{equsw} and \eqref{equsz}, respectively, are both symmetric and unimodal. In particular, for all integers $m,n\ge 0$,
\begin{align*}
&w(m,n)> w(m+2,n),\;\text{provided}\; w(m,n)\neq 0;\\
&z(m,n)> z(m+2,n),\;\text{provided}\; z(m,n)\neq 0.
\end{align*}
\end{theorem}

\section{Unimodality and certain Hilbert scheme}\label{sec5}

In this section, we examine the unimodality of the Betti numbers and Gromov-Witten invariants of certain Hilbert schemes of points related to integer partitions appeared in  G\"{o}ttsche \cite{MR1032930} and Oberdieck \cite{MR3720346}.

\subsection{Betti numbers of certain Hilbert scheme}
In this subsection, following G\"{o}ttsche \cite{MR1032930},
let $X$ be a smooth projective surface over $\bc$. The Hilbert scheme $X^{[n]}:={\sf Hilb}^n(X)$ of points on $X$ is a parameter variety for finite subschemes of length $n$ on $X$. Denote by
\begin{align}\label{eqpoi}
p(X, z)=\sum_{0\le i\le 4} (-1)^ib_i(X)z^i,
\end{align}
the Poinear\'{e} polynomial of $X$, where $b_i(X)$ is the $i$-th Betti number of $X$.

\medskip

In \cite[Theorem 0.1]{MR1032930}, G\"{o}ttsche established the following very nice generating function for the Betti numbers $b_i(X^{[n]})$ of the Hilbert scheme $X^{[n]}$:
\begin{align*}
\sum_{n\ge 0} \sum_{i}b_i(X^{[n]})z^i t^n=\prod_{k\ge 1}\prod_{i=0}^4(1-(-1)^iz^{2k-2+i}t^k)^{(-1)^{i+1}b_i(X)},
\end{align*}
where $b_0(X)=b_{4}(X)$ and $b_1(X)=b_3(X)$. We define that
\begin{equation*}
\mathfrak{G}_{X}(z,q):=\sum_{n\ge 0}\sum_{i\in\bz}b_i(X^{[n]})z^{i-2n} q^n,
\end{equation*}
then the above generating function for $b_i(X^{[n]})$ can be written as:
\begin{align}\label{eqgeq}
\mathfrak{G}_{X}(z,q)=\prod_{k\ge 1}\frac{[(1+zq^k)(1+z^{-1}q^k)]^{b_1}}{(1-q^k)^{b_2}[(1-z^{2}q^k)(1-z^{-2}q^k)]^{b_0}}.
\end{align}
where we write $b_i:=b_i(X)$ for each $i$ for simplicity.

\medskip

Based on Theorem \ref{main1}, the unimodality of $\mathfrak{G}_{X}(z,q)$ will be fully determined by the first factor in the product \eqref{eqgeq}. We can rewrite it as follows:
\begin{align}\label{eqggg}
g_X(z,q):=\frac{[(1+zq)(1+z^{-1}q)]^{b_1}}{(1-q)^{b_2}[(1-z^{2}q)(1-z^{-2}q)]^{b_0}}
=\frac{f_1(z,q)}{(1-q)^{b_2-b_1}}=\frac{f_2(z,q)}{(1-q)^{b_2-b_0}},
\end{align}
where
$$f_1(z,q)=\left(\frac{(1+zq)(1+z^{-1}q)}{(1-q)(1-z^{2}q)(1-z^{-2}q)}\right)^{b_0}\left(\frac{(1+zq)(1+z^{-1}q)}{1-q}\right)^{b_1-b_0},$$
and
$$f_2(z,q)=\frac{\left((1+zq)(1+z^{-1}q)\right)^{b_1}}{[(1-q)(1-z^{2}q)(1-z^{-2}q)]^{b_0}}.$$
To continue our discussion we need the following lemma.
\begin{lemma}\label{lemm12}We have
\begin{equation*}
\frac{(1+q)(1+z q)(1+z^{-1}q)}{(1-z^2q)(1-z^{-2}q)}\in \mathscr{U}_{z,q}^1.
\end{equation*}
\end{lemma}
\begin{proof}
By \eqref{eqmmm0}, we have
$$\frac{1+q}{(1-x^2q)(1-x^{-2}q)}=\sum_{k\ge 0}q^k\sum_{|r|\le k}x^{2r}.$$
Therefore,
\begin{align*}
\frac{(1+q)(1+xq)(1+x^{-1}q)}{(1-x^2q)(1-x^{-2}q)}=&\left(1+q^2+q(x^{-1}+x)\right)\sum_{k\ge 0}q^k\sum_{|r|\le k}x^{2r}\\
=&\sum_{k\ge 0}q^k(1+q^2)\sum_{\substack{|r|\le 2k\\ r\equiv 0\pmod 2}}x^{r}+\sum_{k\ge 1}q^{k}\bigg(\sum_{\substack{|r|\le 2k-2\\ r\equiv 1\pmod 2}}x^{r}+\sum_{\substack{|r|\le 2k\\ r\equiv 1\pmod 2}}x^{r}\bigg)\\
=&\sum_{k\ge 0}q^{k}\left(\sum_{|r|\le 2k}x^{r}+\sum_{|r|\le 2k-3}x^{r}\right)\in \mathscr{U}_{x,q}^1,
\end{align*}
which complete the proof.
\end{proof}
Now, special cases of Lemma \ref{lemm11} imply that
$$\frac{1}{(1-q)(1-z^{2}q)(1-z^{-2}q)}=\frac{1}{(1-q^2)}\frac{1+q}{(1-z^{2}q)(1-z^{-2}q)}\in \mathscr{U}_{z^2,q}^1\subseteq \mathscr{U}_{z,q}^2,$$
and $(1+zq)(1+z^{-1}q)\in \mathscr{U}_{z,q}^2$. Therefore, if $b_2\ge b_0$ then $g_X(z,q)\in \mathscr{U}_{z,q}^2$.
\medskip

On the other hand, a special case of Lemma \ref{lemm11} implies
$$\frac{(1+zq)(1+z^{-1}q)}{1-q}=\frac{(1+q)(1+zq)(1+z^{-1}q)}{1-q^2}\in \mathscr{U}_{z,q}^1,$$
and Lemma \ref{lemm12} implies
$$\frac{(1+zq)(1+z^{-1}q)}{(1-q)(1-z^{2}q)(1-z^{-2}q)}=\frac{1}{1-q^2}\frac{(1+q)(1+zq)(1+z^{-1}q)}{(1-z^{2}q)(1-z^{-2}q)}\in \mathscr{U}_{z,q}^1.$$
Therefore, if $b_2\ge b_1\ge b_0$ then $g_X(z,q)\in \mathscr{U}_{z,q}^1$.

\medskip

Moreover, using Taylor theorem for \eqref{eqgeq} and the definition \eqref{eqpoi} of $p(X,z)$, we have
$$\mathfrak{G}_{X}(z,q)=1+z^{-2}p(X,z)q+O(q^2).$$
Thus, based on the above discussion, Definition \ref{defm0} for $\mathscr{U}_{z,q}^\nu$ and Theorem \ref{main1}, it is not difficult to see the following unimodality of $\mathfrak{G}_{X}(z,q)$.
\begin{theorem}\label{thmsbn}For each $\nu\in\{1,2\}$, $\mathfrak{G}_{X}(z,q)\in \mathscr{U}_{z,q}^\nu$ if and only if  $z^{-2}p(X, z)\in \mathscr{U}_z^\nu$.
\end{theorem}

\begin{remark}
Since $\mathfrak{G}_{X}(z,q)$ is the generating function for the Betti numbers $b_i(X^{[n]})$ of the Hilbert scheme $X^{[n]}$, it is possible that Theorem \ref{thmsbn} can also be proved using the Hard Lefschetz theorem. For details, see Stanley \cite[pp.525-530]{MR1110850}.
\end{remark}

We next use Corollary \ref{prop3210} to study the strict unimodality of $\mathfrak{G}_X(z,q)$, which we present in the following theorem.
\begin{theorem}\label{thmssbn}
For each $\nu\in\{1,2\}$, $\mathfrak{G}_{X}(z,q)\in \mathscr{T}_{z,q}^\nu$ if and only if  $z^{-2}p(X, z)\in \mathscr{T}_z^\nu$.
\end{theorem}
\begin{proof}
Similar to the proof of Theorem \ref{thmsbn}, we only need to show that if $z^{-2}p(X, z)\in \mathscr{T}_z^\nu$ then $g_X(z,q)\in \mathscr{T}_{z,q}^\nu$. In the following, suppose that $z^{-2}p(X, z)\in \mathscr{T}_z^\nu$. To complete the proof of the theorem, it is sufficient to show that
$(1-q)^{-1}f_\nu(z,q)\in \mathscr{T}_{z,q}^\nu$, where $f_{\nu}(z,q)$ is defined as in \eqref{eqggg}. To do this, we write that
$$f_\nu(z,q):=\sum_{n\ge 0}\sum_{i\in \bz}c_{i,n}(\nu)z^iq^n.$$
Notice that for $\nu\in\{1,2\}$, one has
$$f_\nu(z,q)=\frac{\left[(1+zq)(1+z^{-1}q)\right]^{b_1}}{(1-q)^{b_{2-\nu}}[(1-z^{2}q)(1-z^{-2}q)]^{b_0}}\in \mathscr{U}_{z,q}^\nu.$$
Therefore, we shall use Corollary \ref{prop3210} to establish that $(1-q)^{-1}f_{\nu}(z,q)\in \mathscr{T}_{z,q}^\nu$.

\medskip

By applying Newton's binomial theorem to the above expansion for $f_{\nu}(z,q)$ and by using arguments related to Cauchy's products, we have:
\begin{align}\label{eqccc}
c_{i,n}(\nu)=\sum_{\substack{i_0,i_1,i_2,r,s\ge 0\\ i_0+i_1+i_2+r+s=n\\ 2(i_1-i_2)+(s-r)=i}}\binom{b_{2-\nu}-1+i_0}{i_0}\binom{b_0-1+i_1}{i_1}\binom{b_0-1+i_2}{i_2}\binom{b_1}{r}\binom{b_1}{s}.
\end{align}

For the case $b_0=0$, we see that
\begin{align*}
c_{i,n}(\nu)=\sum_{\substack{i_0,r,s\ge 0\\ i_0+r+s=n,\; s-r=i}}\binom{b_{2-\nu}-1+i_0}{i_0}\binom{b_1}{r}\binom{b_1}{s}=\sum_{\substack{i_0,r\ge 0\\ i_0+2r=n-i}}\binom{b_{2-\nu}-1+i_0}{i_0}\binom{b_1}{r}\binom{b_1}{r+i}.
\end{align*}
Therefore, $c_{i,n}(\nu)=0$ for all $i>\min(n,b_1)$ and
\begin{align*}
c_{i,i}(\nu)=\binom{b_{2-\nu}-1}{0}\binom{b_1}{0}\binom{b_1}{i}=\binom{b_1}{i}>0=c_{i+\nu,i}(\nu),
\end{align*}
for all $0\le i\le b_1$.  Thus the use of Corollary \ref{prop3210} immediately yields the proof of the theorem for the case $b_0=0$.

\medskip

For the cases $b_0\ge 1$, note that for all integers $t\le 1$ and $i\ge 0$ such that $2i-t\ge 0$,
$$\sum_{\substack{i_0,i_1,i_2,r,s\ge 0\\ i_0+i_1+i_2+r+s=i\\ 2(i_1-i_2)+(s-r)=2i-t}}=\sum_{\substack{i_0,i_1,i_2,r,s\ge 0\\ i_0+i_1+i_2+r+s=i\\ 2i_0+4i_2+3r+s=t}}=\sum_{\substack{s\ge 0\\ i_0=i_2=r=0\\ s=t, i_1=i-\mu}},
$$
and applying to \eqref{eqccc} yields
\begin{align*}
c_{2i-t,i}(\nu)=\binom{b_0-1+i-t}{i-t}\binom{b_1}{t}>0,
\end{align*}
for $t=\{0,1\}$, and $c_{2i-t,i}(\nu)=0$ for all $t<0$.  This implies that  $c_{i,n}(\nu)>0$ for all $0\le i\le 2n$ and $c_{i,n}(\nu)=0$ for all $i>2n$. Therefore,
\begin{align*}
c_{2i-t,i}(\nu)-c_{2i-t+\nu,i}(\nu)=c_{2i-t,i}(\nu)>0,
\end{align*}
for $\nu=2$ with $t\in\{0,1\}$, and $\nu=1$ with $t=0$. For $\nu=1$ with $t=1$, notice that $2i-1\ge 0$, we have $i\ge 1$ and $b_1>b_0$, thus
\begin{align*}
c_{2i-1,i}(1)-c_{2i,i}(1)=\binom{b_0-1+i-1}{i-1}b_1-\binom{b_0-1+i}{i}=\frac{i(b_1-1)-(b_0-1)}{b_0-1+i}\binom{b_0-1+i}{i}>0.
\end{align*}
Combining the above arguments, since we can write $h=2i-t$ with $t\in\{0,1\}$ for any integer $0\le h\le 2n$, we have
$$c_{h,\lceil h/2\rceil}(\nu)-c_{h+\nu,\lceil h/2\rceil}(\nu)=c_{2i-\mu,i}(\nu)-c_{2i-\mu+\nu,i}(\nu)>0.$$
Therefore, the use of Corollary \ref{prop3210} immediately yields the proof of the theorem for the cases $b_0\ge 1$. This completes the proof of the theorem.
\end{proof}

\subsection{Gromov-Witten invariants of certain Hilbert scheme}In this subsection, let $y= - e^{2 \pi i z}$, $q = e^{2 \pi i \tau}$ with $z\in \bc$ and $\Im(\tau)>0$, and let
\begin{equation*}
\eta(\tau) = q^{1/24}(q;q)_\infty
\end{equation*}
denote the Dedekind eta function, and
\begin{equation*}
\vartheta_1(z,\tau) = q^{1/8}(y^{-1/2}+y^{1/2})(q,-yq,-y^{-1}q;q)_\infty
\end{equation*}
denote the first Jacobi theta function.

\medskip

In the recent work \cite{MR3720346}, Oberdieck calculated genus $0$ Gromov-Witten invariant of the Hilbert scheme ${\sf Hilb}^d(S)$ of $d$ points of an elliptically fibered $K3$ surface $ S\to \bp^1$, which count rational curves incident to two generic fibers of the induced Lagrangian fibration ${\sf Hilb}^d(S)\to \bp^d$. Oberdieck \cite[Theorem 1]{MR3720346} proved that the generating series of these invariant is the Fourier expansion of a \emph{Jacobi form}. The generating series in \cite[Theorem 1]{MR3720346} is
\begin{align*}
\sum_{h\ge 0}\sum_{k\in\bz} \mathsf{N}_{d,h,k}y^kq^{h-1}=\eta(\tau)^{-24}\left(\frac{\vartheta_1(z,\tau)}{\eta(\tau)^3}\right)^{2d-2},
\end{align*}
where the Gromov-Witten invariants $\mathsf{N}_{d,h,k}$ count the number of all rational curves in a class depending on the perimeters $h$ and $k$, see Oberdieck \cite{MR3720346} for details.

\medskip

Oberdieck \cite{MR3720346} also calculated genus $0$ Gromov-Witten invariants for several other natural incidence conditions.
He \cite[Theorem 2]{MR3720346} proved that the generating series of these invariant are also Jacobi forms. More precisely, define that
$$G(z,\tau):=\eta(\tau)^{-6}\vartheta_1(z,\tau)^2\left(y \frac{\rd}{\rd y}\right)^2 \log \vartheta_1(z,\tau).$$
Then the generating series in \cite[Theorem 2]{MR3720346} are
\begin{align*}
\sum_{h\ge 0}\sum_{k\in\bz} \mathsf{N}_{d,h,k}^{(1)}y^kq^{h-1}&=\eta(\tau)^{-24}G(z,\tau)^{d-1},\\
\sum_{h\ge 0}\sum_{k\in\bz} \mathsf{N}_{d,h,k}^{(2)}y^kq^{h-1}&=\frac{1}{2(1-d)}\eta(\tau)^{-24}\left(y \frac{\rd}{\rd y} G(z,\tau)^{d-1}\right),\nonumber\\
\sum_{h\ge 0}\sum_{k\in\bz} \mathsf{N}_{d,h,k}^{(3)}y^kq^{h-1}&=\frac{1}{d}\binom{2d-2}{d-1}\eta(\tau)^{-24}\left(q \frac{\rd}{\rd q}G(z,\tau)^{d-1}\right),
\end{align*}
where the Gromov-Witten invariants $\mathsf{N}_{d,h,k}^{(i)}, (i=1,2,3)$ count the number of all rational curves in a class depending on the perimeters $i$, $h$ and $k$, see also Oberdieck \cite{MR3720346} for details.

\medskip

Based on our previous results, we establish the following unimodality for the Gromov-Witten invariants $\mathsf{N}_{d,h,k}$, $\mathsf{N}_{d,h,k}^{(1)}$ and $\mathsf{N}_{d,h,k}^{(3)}$.
\begin{theorem}\label{thm311}
For any integers $d\ge 2$ and $h\ge 0$, the sequence $\left(\mathsf{N}_{d,h,k}\right)_{k\in\bz}$ is symmetric and strictly unimodal. For any integers $d\ge 2$ and $h\ge 0$, the sequences $\left(\mathsf{N}_{d,h,k}^{(1)}\right)_{k\in\bz}$ and $\left(\mathsf{N}_{d,h,k}^{(3)}\right)_{k\in\bz}$ are symmetric and unimodal.
\end{theorem}
It is not difficult to see that Theorem \ref{thm311} follows from Theorem \ref{main1} and the following Proposition    \ref{pro312}.
\begin{proposition}\label{pro312}We have
$$q^{-1/24}\eta(\tau)^{-5}\vartheta_1(z,\tau)^2\in \mathscr{T}_{y,q}^1,$$
and
$$q^{-1/4}\eta(\tau)^{-4}\eta(2\tau)^2\vartheta_1(z,\tau)^2\left(y \frac{\rd}{\rd y}\right)^2 \log \vartheta_1(z,\tau)\in \mathscr{U}_{y,q}^1.$$
\end{proposition}
\begin{proof}Clearly, the fact
$q^{-1/24}\eta(\tau)^{-5}\vartheta_1(z,\tau)^2\in \mathscr{T}_{y,q}^1$ is a special case of Theorem \ref{thmssbn}. It remains to prove the second relation. Note that
\begin{align*}
\left(y \frac{\rd}{\rd y}\right)^2 \log \vartheta_1(z,\tau) &=\frac{y}{(1+y)^2}+\sum_{m\ge 1}\left(\frac{y^{-1}q^m}{(1+y^{-1} q^m)^2}+\frac{yq^m}{(1+y q^m)^2}\right)\\
&=\frac{y}{(1+y)^2}+\sum_{m\ge 1}\frac{\left(4q^{m}+(y^{-1}+y)(1+q^{2m})\right)q^m}{(1+y^{-1} q^m)^2(1+y q^m)^2},
\end{align*}
and
\begin{align*}
\frac{(1+y)^2}{y}\left(4q+(y^{-1}+y)(1+q^{2})\right)=\left(\frac{1}{y^2}+\frac{2}{y}+2+2y+y^2\right)(1+q^2)+4\left(\frac{1}{y}+2+y\right)q\in \mathscr{U}_{y,q}^1,
\end{align*}
we have
\begin{align*}
&\frac{\eta(2\tau)^2\vartheta_1(z,\tau)^2}{q^{1/4}\eta(\tau)^{4}}\left(y \frac{\rd}{\rd y}\right)^2 \log \vartheta_1(z,\tau)\\
=&(-q,-yq,-y^{-1}q;q)_\infty^2+\sum_{m\ge 1}\frac{(1+y)^2}{y}\left(4q^m+(y^{-1}+y)(1+q^{2m})\right)\frac{q^m(-q,-yq,-y^{-1}q;q)_\infty^2}{(1+y^{-1} q^m)^2(1+y q^m)^2}.
\end{align*}
Thus we immediately obtain the proof of this proposition by using Theorems \ref{main1} and \ref{main}.
\end{proof}

\section*{Acknowledgements}
The author would like to thank the anonymous referees for their very helpful comments and suggestions.


\end{document}